\documentclass[12pt,a4paper]{amsart}
\usepackage[left=2.85cm,right=2.85cm,top=2.54cm,bottom=2.54cm]{geometry}
\usepackage{amsmath, amssymb, amsthm, amsbsy, mathrsfs}
\usepackage[colorlinks=true, linkcolor=blue, citecolor=blue, urlcolor=blue]{hyperref}
\usepackage{microtype}
\usepackage{tikz}

\newtheorem{theorem}{Theorem}

\theoremstyle{definition}

\newcommand{\F}{\mathfrak F}
\newcommand{\other}{\beta}
\renewcommand{\P}{\mathbb P}
\newcommand{\sameF}{\alpha}
\newcommand{\sameP}{\gamma}

\begin{document}
\title[Friedman vs P\'olya]{Friedman vs. P\'olya}
\date{\today}
\author[R. Alves]{Raphael Alves}
\author[R. A. Rosales]{Rafael A. Rosales}
\address[R. Alves, R. A. Rosales]{%
  Departamento de Computa\c{c}\~ao e Matem\'atica,
  Universidade de S\~ao Paulo,
  Ribeir\~ao Preto, 14040-901 - SP, Brasil
}
\email{raphael.alves.duarte@usp.br}
\email{rrosales@usp.br}
\thanks{The first author was partially supported by 
  a scholarship from the \emph{Programa Unificado de Bolsas da
    Universidade de S\~ao  Paulo}, \texttt{\#}2023-4061.
} 

\begin{abstract}
  Suppose an urn contains initially any number of balls of two
  colours. One ball is drawn randomly and then put back with $\alpha$
  balls of the same colour and $\beta$ balls of the opposite
  colour. Both cases, $\beta=0$ and $\beta>0$ are well known and
  correspond respectively to P\'olya's and Friedman's replacement
  schemes. We consider a mixture of both of these: with probability
  $p\in(0,1]$ balls are replaced according to Friedman's recipe and with
  probability $1-p$ according to the one by P\'olya. Independently of
  the initial urn composition and independently of $\alpha$, $\beta$,
  and the value of $p>0$, we show that the proportion of balls of one
  colour converges almost surely to $\frac12$. The latter is the limit
  behaviour obtained by using Friedman's scheme alone, i.e. when
  $p=1$. Our result follows by adapting an argument due to
  D. S. Ornstein.
\end{abstract}



\begingroup
 \def\uppercasenonmath#1{\scshape} 
 \let\MakeUppercase\relax 
 \maketitle
 \endgroup


\section{Introduction}
An urn contains $Y_n$ yellow balls and $B_n$ blue balls at time
$n$. One ball is drawn randomly and then put back together with
$\sameF$ balls of the same colour and $\other$ balls of the opposite
colour. The case $\other=0$ is the well known P\'olya urn model in
which $X_n = Y_n/(Y_n+B_n)$ converges a.s. to a limiting random variable
that has a beta distribution with parameters $Y_0/\sameF$ and
$B_0/\sameF$. On the other hand, the case $\other>0$ behaves radically
differently: independently on the initial composition $Y_0$, $B_0$ and
no matter how large $\sameF$ is in comparison with $\other$,
$X_n$ converges always a.s. towards $\frac12$. This case has
also been extensively studied and is known as Friedman's urn or
equivalently as a generalised P\'olya urn.

Suppose now that at each step, with probability $p$ one chooses
Friedman's replacement scheme and with probability $1-p$ the one by
P\'olya. In this case, we show that for any initial configuration
$Y_0$ and $B_0$ and any positive $p$, no matter how small, the limit
of $X_n$ follows the one obtained by using Friedman's scheme
alone. The histograms for $X_n$ shown below, are obtained from three
simulations made by taking $B_0=Y_0=1$ and
$\alpha=\beta=1$.\\[-.5em]

\begin{center} 
 \input{f_vs_p-V2} 
\end{center}
The plot at the left shows the result obtained with $p=0$. The ones at
the center and the right were made with $p=0.05$ for $2\times 10^3$
and $2\times10^7$ draws respectively. They suggest that $X_n$
concentrates at $\frac12$ when $p=0.05$, although very slowly.

Our main result, stated below, is slightly more general than the one
described in the previous paragraph because we allow the addition of
balls of the observed colour via Friedman's or P\'olya's scheme to be
different. Hereafter, we denote these numbers by $\sameF$ and $\sameP$
respectively.

\begin{theorem}\label{th:friedman_wins}
For any positive integers $Y_0$, $B_0$, $\sameF$, $\other$, $\sameP$,
and any $p \in (0, 1]$, almost surely $X_n \longrightarrow\frac12$. 
\end{theorem}

The rest of this note is devoted to the proof
Theorem~\ref{th:friedman_wins}.  The main argument proceeds by
adapting the proof of the convergence of Friedman's urn as presented
in \cite[Example 5.4.5]{D10}. The general idea behind this can be
traced back to Donald S. Ornstein, see \cite{Freedman65}. The theorem
may also be established by other means such as stochastic
approximations or the analytic approach to generalized P\'olya
urns, see \cite{LP19} or \cite{FDP06}, \cite{MM12}. We are however not
aware of the existence of any such proofs.

\section{Proof of Theorem~\ref{th:friedman_wins}}
Let $(Y_n, B_n)$, $n\geq 0$, be defined on a probability
space $(\Omega, \F, \P)$ endowed with $\F_n$, the natural filtration
generated by $(Y_n, B_n)$. Let $(\zeta)_{n\geq 1}$ be a sequence of
random variables such that $\P(\zeta_n = 1\mid \F_{n-1}) = 1 -
\P(\zeta_n = 0\mid \F_{n-1}) = X_{n-1}$.  Let $(\xi_n)_{n\geq 1}$ be
another sequence of random variables, independent of $\F_n$ and
$(\zeta_n)_{n\geq 1}$, defined such that $\P(\xi_n=1) = 1 -
\P(\xi_n = 0) = p$, $p \in (0,1]$. In this case
\begin{align*}
  Y_n
  &= Y_{n-1} + \sameF \mathbf{1}_{\{\xi_n=1,\ \zeta_n=1\}} +
  \other\mathbf{1}_{\{\xi_n=1,\ \zeta_n=0\}} +
  \sameP\mathbf{1}_{\{\xi_n=0,\ \zeta_n=1\}} \\
  &= Y_0 + \alpha \sum_{k=1}^n \mathbf{1}_{\{\xi_k=1,\ \zeta_k=1\}} 
    + \other\sum_{k=1}^n \mathbf{1}_{\{\xi_k=1,\ \zeta_k=0\}} +
    \sameP\sum_{k=1}^n \mathbf{1}_{\{\xi_k=0,\ \zeta_k=1\}}.
\end{align*}

We will make use of the following version of the Second Borel-Cantelli
Lemma, see \cite[Theorem 5.4.11]{D10}. Suppose $A_n$ is a sequence of
events adapted to $\F_n$ and let $p_n = \P(A_n \mid \F_{n-1})$. Then
\[
  \sum_{k=1}^n \mathbf{1}_{A_k} \Bigg/ \sum_{k=1}^n p_{k-1}
  \longrightarrow
  1\
  \quad\text{a.s. on}\quad
  \Bigg\{\sum_{k=1}^\infty p_{k-1} = \infty\Bigg\}.
\]

Let $A_{ijk} = \{\xi_k = i, \zeta_k = j\}$ and $\sigma =
\max\{\alpha+\beta, \gamma\}$. It then follows that for any $k$, $p_k
= Y_k/(Y_k+B_k) \geq Y_0/(Y_0+B_0+\sigma k)$. Thus, for $(i,j)$ equal
to either $(1,1)$, $(0,1)$ or $(1,0)$, we have that $\sum_{k=1}^\infty
\P(A_{ijk}\mid\F_{k-1}) = \infty$ and so, the conditions of the
Borel-Cantelli Lemma are satisfied by any one of the sequences
$(A_{ijn})_{n\geq 1}$. Then
\begin{align*}
  \limsup_n \frac{Y_n}{n}
  &=
  \limsup_n\frac1n
  \sum_{k=1}^n \Big\{
   \sameF\, \P\big(\xi_k = 1, \zeta_k=1 \mid \F_{k-1}\big)
    + \other\, \P\big(\xi_k = 1, \zeta_k = 0 \mid \F_{k-1}\big) \\ 
  &\qquad\qquad\qquad\qquad +
   \sameP \P\big(\xi_k = 0, \zeta_k = 1\mid \F_{k-1}\big)\Big
  \},
\end{align*}
which, by the independence of $\xi_k$ and $\zeta_k$ gives
\begin{align*}
    \limsup_n \frac{Y_n}{n}
  &=
    \limsup_n \frac1n \sum_{k=1}^n \bigg\{ \alpha
    p\, \P\big(\zeta_k=1\mid\F_{k-1}\big) 
    + \beta p\,  \P\big(\zeta_k=0\mid\F_{k-1}\big)\\
    &\qquad\qquad\qquad\qquad + \gamma(1-p)
      \P\big(\zeta_k=1\mid\F_{k-1}\big) 
      \bigg\} \\
  &=
    \limsup_n \frac1n \sum_{k=1}^n \bigg\{ \left[\alpha p-\beta
    p+\gamma(1-p)\right] X_{k-1} 
    + \beta p\bigg\}.
\end{align*}
Similarly
\begin{align*}
  \limsup_n \frac{Y_n+B_n}{n}
  &=
    \limsup_n \frac1n \sum_k\bigg\{
    (\alpha+\beta) \P\big(\xi_k = 1\mid\F_{k-1}\big)+
    \gamma\, \P(\xi_k=0\mid\F_{k-1})\bigg\} \\
  &=
    (\alpha+\beta)p + \gamma(1-p).
\end{align*}
The previous computations lead finaly to
\begin{equation}
  \label{eqn:Xn_recursion}
  \limsup_n X_n = \limsup_n \frac{\big[(\alpha - \beta) p +
    \gamma(1-p)\big]X_{n-1} + \beta p}{(\alpha + \beta)p +
    \gamma(1-p)}.
\end{equation}

Let $\theta = (\alpha-\beta)p + \gamma(1-p)$, and let $\ell: [0,1] \to
[0,1]$ be defined by
\begin{equation}
  \label{eqn:reta}
  \ell(x) = \frac{\theta x + \beta p}{(\alpha+\beta)p + \gamma(1-p)}.
\end{equation}
For any $\alpha, \beta, \gamma$ and $p$, $\ell$ is a linear function
in $x$, with slope taking values in $(-1,1)$ and unique fixed point at
$\frac12$. The slope of $\ell$ may be positive or negative according
to the sign of $\theta$. We consider the
possible cases separately, inclusive when the slope is zero.

\textbf{Case 1}: $\theta > 0$. Observe that $\ell$ is increasing with
$x$. Assume that $\limsup_n X_{n-1} \leq x$, for $x \in [0,1]$. From
\eqref{eqn:reta} this gives the upper bound $\limsup_n X_n \leq
\ell(x)$.  Taking $x=1$ as an obvious initial guess for $x$, and
substituting this into the right hand side of the previous inequality
shows that $\limsup_{n} X_n \leq \ell(1) < 1$. Iterating leads to
$\limsup_n X_n \leq \lim_n \ell^{(n)}(1) = \frac12$. A similar
argument obtained by interchanging colours gives $\liminf_n X_n \geq
\frac12$ and concludes the proof.

\textbf{Case 2}: $\theta < 0$. In this case,
the right hand side of \eqref{eqn:reta} is maximised when
$x=0$ and the argument used in the previous case does not go
through. A way out is obtained by considering the proportion of blue
balls $Z_n$ together with $X_n$, the proportion of yellow
balls. Since $Z_n=1-X_n$, straightforward computations starting from
\eqref{eqn:Xn_recursion} give 
\begin{align*}
  \limsup_n Z_n
  &=
    \limsup_n \frac{-\theta X_{n-1} + \alpha p +
    \gamma(1-p)}{(\alpha+\beta)p+\gamma(1-p)}, \\
  \limsup_n X_n
  &=
    \limsup_n \frac{-\theta Z_{n-1} + \alpha p +
    \gamma(1-p)}{(\alpha+\beta)p+\gamma(1-p)},
\end{align*}
Assume now that $\limsup_n Z_{n-1}\leq z$ and $\limsup_n X_{n-1} \leq
x$, for some $z, x \in [0,1]$, which then leads to
\[
  \limsup_n Z_n
  \leq
    \frac{-\theta x + \alpha p +
  \gamma(1-p)}{(\alpha+\beta)p+\gamma(1-p)}
  \quad\text{and}\quad
  \limsup_n X_n
  \leq
    \frac{-\theta z + \alpha p +
    \gamma(1-p)}{(\alpha+\beta)p+\gamma(1-p)}.
\]
Starting with the initial values $x=1$ and $z=1$, and observing that
$-\theta > 0$ allows one to proceed exactly as in Case 1 to conclude
that both $\limsup_n Z_n$ and $\limsup_n X_n$ are less or equal to
$\frac12$. Both $\liminf_n Z_n$ and $\liminf_n X_n$ are also easily
shown to be greater or equal to $\frac12$.

\textbf{Case 3}: $\theta = 0$. Isolating $p$ in this case leads to $p
= \gamma/(\beta+\gamma-\alpha)$. Substituting this expression into
\eqref{eqn:reta} directly establishes $\limsup_n X_n \leq
\frac12$. Proceeding in the same fashion but changing colours shows
that $\liminf_n X_n \geq \frac12$. This concludes the proof of the
theorem.


\vspace{.75cm} 

\end{document}